\documentclass[11pt,fleqn]{scrartcl} 
  
\usepackage[utf8]{inputenc}
\usepackage[T1]{fontenc}

\usepackage{lmodern}
\usepackage{amsmath,amssymb}
\usepackage{graphicx}
\usepackage{dsfont}
\usepackage{exscale}
\usepackage{hyperref}
\usepackage{tikz-cd}
\usepackage{float}
\usepackage[lofdepth,lotdepth]{subfig}
\usepackage[ngerman]{translator}
\usepackage{listings}

\usepackage[backend=biber, %% Hilfsprogramm "biber" (statt "biblatex" oder "bibtex")
sorting=nyt, 
style=alphabetic, %% Zitierstil (siehe Dokumentation)
hyperref=true, %% hyperref-Paket verwenden, um Links zu erstellen
isbn=false, %keine isbn zeigen
url=true,
]{biblatex}
 
\newcommand{\eop}{\hfill $\square$}
\newtheorem{defi}{Definition}[section]
\newtheorem{satz}[defi]{Theorem}

\newcommand{\D}{\mathbf{\mathcal{D}}}
\newcommand{\X}{\mathbf{x}}
\newcommand{\uX}{\underline{\X}}

\newtheorem{example}{Example}

\begin{document}
\title{First steps towards $q$-deformed Clifford analysis}
\author{M. Zimmermann\footnote{TU Bergakademie Freiberg, Institut of Applied Analysis}, S. Bernstein$^*$, B. Schneider\footnote{University of  Ostrava, Department of Mathematics} }
\date{}
\maketitle
\tableofcontents
\begin{abstract} Abstract. We consider the extension of the Jackson calculus into higher dimensions and specifically into Clifford analysis. \end{abstract}
\section{Introduction}
In general, q-deformed stands for quantum-deformed. Essentially, one can see there different directions here. One is the $q$-deformed space and the $q$-deformed sphere, the fuzzy sphere. The q-deformed spaces, in turn, lead to non-commutative structures and hence to the q-analogues of known Lie groups. There are interesting applications in quantum physics and optics. 
However, this article is not about physics or physical mathematics, but about the closely related Jackson calculus. Jackson calculus replaces the ordinary derivative $$f^{\prime}(x_0) = \lim_{x\to x_0}\frac{f(x)-f(x_0)}{x-x_0}$$ with the difference operator $$D_qf(x) = \frac{d_qf(x)}{d_qx}= \frac{f(qx)-f(x)}{qx-x}, $$ where $\lim_{q\to 1}D_qf(x) = f^{\prime}(x). $ Jackson calculus is described in the book by Jackson \cite{Jackson1909} or in the more recent book \cite{KacCheung2002}. References to quantum physics and Jackson calculus can be found in \cite{ChaiDem1996}, \cite{Ernst2012}. 
The present paper is concerned with the extension of the Jackson calculus into higher dimensions and specifically into Clifford analysis \cite{Brackx1982, Delanghe2012, GS1989, GHS2016, GM1991}. There is already an extension of the Jackson calculus into Clifford analysis by Coulembier and Sommen \cite{Coulembier2010, Coulembier2011}. In these papers, a Dirac operator is defined by axioms. In the present work, a Jackson calculus or q-calculus is introduced by defining the partial derivatives by $$
\partial_{x_i}^q(f(x_1,...,x_n)) = \frac{f(x_1,...,x_{i-1},qx_i,x_{i+1},...x_n) - f(x_1,...,x_n)}{(q-1)x_i}$$ and thus $q$-forms. This is then the basis for the $q$-deformed Dirac operator $\D_{\uX}^q = - \sum_{i=1}^m e_i \partial_{x_i}^q.$ This Dirac operator does not satisfy the axioms of the Dirac operator of Coulembier and Sommen. Therefore, in this short article, we first consider the $q$-defined Dirac, Euler and Gamma operators and their symmetry properties. Then the Fisher decomposition and the corresponding Cauchy-Kovalevskya expansion with examples.  

\section{Preliminaries}
\subsection{Clifford analysis}
We consider the Clifford algebra $\mathcal{C}\ell_{0,m}$ over $\mathbb{R}^m$ with the identity element $e_{\emptyset}$ satisfying $e_{\emptyset}^2 = 1$ and elements $e_1,...,e_m$ following the multiplication rules
\begin{equation}
e_ie_j + e_je_i = -2\delta_{ij}
\end{equation}
for $i,j = 1,...,m$. The $e_i$ form a basis with $2^m$ elements: $(e_{\emptyset},e_1,e_2,...,e_1e_2,...,e_1e_2e_3,...)$. Using $M:= \{1,...,m\}$ and $A:= \{(h_1,...,h_r \in \mathcal{P}M : 1 \leq h_1 \leq ... \leq h_r \leq m \}.$ An arbitrary element of the Clifford algebra $\lambda \in \mathcal{C}\ell_{0,m}$ is given by
\begin{equation}
\lambda = \sum_A \lambda_A e_A\text{, } \lambda_A \in \mathbb{R}
\end{equation}
with $e_A = e_{h_1}...e_{h_r}$. We consider Clifford-valued functions $f$, i.e.
\begin{equation}
f(\uX) =\sum_A f_A(\uX)e_A = \sum_A f_A(x_1,...,x_m)e_A\text{, } f_A \colon \mathbb{R}^m \rightarrow \mathbb{R}.
\end{equation}
The vector variable is defined as
\begin{equation}
\uX = \sum_{i=1}^m x_i e_i.
\end{equation}
The multiplication rules imply that $\uX^2$ is scalar valued, $\uX^2 = -\sum_{i=1}^m x_i^2 = - \vert \uX \vert^2$. Further we can now define the Dirac operator
\begin{equation}
\D_{\uX} = - \sum_{i=1}^m e_i \partial_{x_i},
\end{equation}
where $(\D_{\uX})^2 = - \Delta$ is scalar and $\Delta$ denotes the Laplace operator.
A $\mathcal{C}\ell_{0,m}$-valued function $f = \sum_A e_Af_A(\uX)$ is called \textbf{left (right) monogenic} if and only if
\begin{equation}
\D_{\uX}f = 0\quad (f\D_{\uX} = 0).
\end{equation}
More on Clifford analysis can be found in the following monographs \cite{Brackx1982}, \cite{Delanghe2012}.

\subsection{Spherical harmonics and monogenics}
%If $\D_{\uX}f = f\D_{\uX} = 0$ applies, the function is called monogenic. 
%\textbf{Polynomials of degree at most $k$}\\
The space of \textbf{polynomials of degree at most $k$} is
\begin{equation}
\Pi_k = \left\{ P(\uX): P(\uX) = \sum_{\vert \alpha \vert \leq k} c_{\alpha}\uX^{\alpha} \right\}
\end{equation}
for $\uX \in \mathbb{R}^m, \, \uX = (x_1,...,x_m)$ and an index vector $\alpha = [\alpha_1,...,\alpha_m] \in \mathbb{N}^m$ with $\vert \alpha \vert = \alpha_1 + ... +\alpha_m$ and $x^{\alpha} = x_1^{\alpha_1}x_2^{\alpha_2}...x_m^{\alpha_m}$.
\\
A polynomial $P(\uX)$ is called homogeneous of degree $k$ if it satisfies for every $\uX \in \mathbb{R}^m$ and $\uX \neq 0$ that
\begin{equation}
P(\uX) = \vert \uX \vert^k P\left(\frac{\uX}{\vert \uX \vert}\right).
\end{equation}
A homogeneous polynomial of degree $k$ is an eigenfunction of the Euler operator $\mathbb{E}$ with eigenvalue $k$. Using this we can define the space of homogeneous polynomials of degree $k$  
\begin{equation}
\mathcal{P}_k = \{ P(\uX): \mathbb{E}P(\uX) = kP(\uX)\}.
\end{equation}
%\textbf{Harmonics and spherical harmonics}\\
Harmonic functions are the null-solutions of the Laplace operator $\Delta = \sum_{i=1}^m \partial_{x_i}^2$. Therefore the space of $k$-homogeneous polynomials that are also harmonic is denoted by
\begin{equation}
\mathcal{H}_k = \{P(\uX) \in \mathcal{P}_k : \Delta P(\uX) = 0\}.
\end{equation}
The restriction of a polynomial $P(\uX) \in \mathcal{H}_k$ to the unit sphere $\mathbb{S}^{m-1}$ is called a spherical harmonic of degree $k$. \\
\\
%\textbf{$\mathcal{C}\ell_m$-valued homogeneous polynomials}\\
A $\mathcal{C}\ell_{0,m}$-valued homogeneous polynomial is a polynomial
\begin{align}
P(\uX) &= \sum_A c_A\uX^A,\, c_A \in \mathcal{C}\ell_{0,m} \\
\Leftrightarrow P(\uX) &= \sum_A e_Ap_A(\uX)
\end{align}
with $p_A(\uX)$ a complex-valued $k$-homogeneous polynomial and $\vert A \mid = k$. \\
\\
%\textbf{Monogenic, homogeneous polynomials}\\
The space of monogenic homogeneous polynomials of degree $k$ is denoted by
\begin{equation}
\mathcal{M}_k = \{ M(x) \in \mathcal{P}_k : \D_{\uX}M = 0\}.
\end{equation}
That is monogenic polynomials that are also eigenfunctions of the Euler operator. The restriction of a polynomial $M(x) \in \mathcal{M}_k$ to the unit sphere is called a spherical monogenic of degree $k$. \\

We want to extend Jackson's calculus to higher dimensions. Therefore we recall some of the one-dimensional Jackson calculus.
\subsection{Jackson calculus}
For Jackson's calculus see also \cite{Jackson1909, Ernst2012}.
For a number $u$ and a deformation parameter $q \in \mathbb{R}^+$ we can define the $q$-deformation of $u$ as
\begin{equation}
[u]_q = \frac{q^u-1}{q-1} = 1+q+q^2+...+q^{u-1}.
\end{equation}
Taking the limit this satisifies $\lim_{q\rightarrow 1} [u]_q = u$. The $q$-derivative or Jackson derivative of a function $f(t)$ is defined as
\begin{equation}
\partial_t^q(f(t)) = \frac{f(qt)-f(t)}{(q-1)t}.
\end{equation}
This yields
\begin{equation}
\partial_t^q(t^k) = [k]_q t^{k-1}.
\end{equation}
Further the Leibniz rule is satisfied
\begin{equation}
\partial_t^qt = qt\partial_t^q+1 
\end{equation}
or more in general
\begin{align}
\partial_t^q(f_1(t)f_2(t)) &= \partial_t^q(f_1(t))f_2(t) + f_1(qt)\partial_t^q(f_2(t))\\
&= \partial_t^q(f_1(t))f_2(qt) + f_1(t)\partial_t^q(f_2(t)),
\end{align}
i.e. there are two versions of a Leibniz rule. \\
We can also define the $q$-integration. For $q<1$ the $q$-integral is given by
\begin{equation}
\int_0^a f(x) d_qx = (1-q)a\sum_{k=0}^{\infty} f(aq^k)q^k.
\end{equation}
Integration on general intervals $[a,b]$ is defined by $\int_a^b = \int_0^b - \int_0^a$. The $q$-binomial coefficient is
\begin{equation}
\binom{n}{k}_q = \frac{[n]_q!}{[n-k]_q![k]_q!}
\end{equation}
with the $q$-factorial $[k]_q! = [k]_q[k-1]_q...[1]_q$ and $[0]_q := 1$. Now we can also define the $q$-exponential
\begin{equation}
E_q(t) = \sum_{j=0}^\infty \frac{t^j}{[j]_q!}
\end{equation}
and for the inverse we need a second $q$-exponential given by
\begin{equation}
e_q(t) = E_{q^{-1}}(t).
\end{equation}
These satisfy the following relations:
\begin{gather}
E_q(t)e_q(-t) = 1\\
\partial_t^qE_q(t) = E_q(t)\\
\partial_t^qe_q(t) = e_q(qt).
\end{gather}

\subsection{$q$-deformed partial derivatives}
As we will consider functions dependent on multiple variables we need some sort of $q$-partial derivative. Following Demichev \cite{Demichev1996} these will be defined by
\begin{equation}
\partial_{x_i}^q(f(x_1,...,x_n)) = \frac{f(x_1,...,x_{i-1},qx_i,x_{i+1},...x_n) - f(x_1,...,x_n)}{(q-1)x_i}.
\end{equation}
The $q$-partial derivative satisfies the following relations:
\begin{itemize}
\item $\partial_{x_j}^q x_i = x_i \partial_{x_j}^q$ for $i\neq j$
\item $\partial_{x_i}^qx_i = qx_i\partial_{x_i}^q+1$
\item $\partial_{x_i}^q \partial_{x_j}^q = \partial_{x_j}^q \partial_{x_i}^q$ for $i \neq j$
\item $\partial_{x_i}^q x_i^2 = q^2x_i^2\partial_{x_i}^q + [2]_qx_i$
\item $(\partial_{x_i}^q)^2x_i = q^2x_i^2(\partial_{x_i}^q)^2 + [2]_q \partial_{x_i}^q$
\item $(\partial_{x_i}^q)^2x_i^2 = q^4x_i^2 (\partial_{x_i}^q)^2 + (q^2+1)[2]_qx_i\partial_{x_i}^q + [2]_q$
\end{itemize}

Now, we can consider an appropriate Clifford analysis by using the $q$-deformed partial derivatives.
\section{$q$-deformed Clifford analysis}
\subsection{$q$-deformed Dirac operator}
First we go back to the Dirac operator defined earlier and replace the continuous partial derivates by the $q$-deformed partial derivatives. This results in a $q$-Dirac operator
\begin{equation}
\D_{\uX}^q = - \sum_{i=1}^m e_i \partial_{x_i}^q.
\end{equation}
A generalization of the q-deformation and the Jackson calculus was already considered in Coul./Sommen. There is an approach based on four axioms which are the basis for their calculation.
According to Coulembier and Sommen \cite{Coulembier2011, Coulembier2010} following axioms should be satisfied
\begin{itemize}
\item (A1) $\D_{\uX}^q(\uX) = [m]_q$
\item (A2) $\D_{\uX}^q\uX^2 = q^2\uX^2\D_{\uX}^q+(q+1)\uX$
\item (A3) $(\D_{\uX}^q)^2$ scalar
\item (A4) $\D_{\uX}^qM_k = 0$.
\end{itemize}
Our approach is different. By defining the $q$-Dirac operator as seen above (A1) and (A2) are not satisfied. Instead we have for (A1)
\begin{equation}
\D_{\uX}^q(\uX) = [1]_q\sum_{i=1}^m 1 = m
\end{equation}
and (A2) has to be replaced by
\begin{equation}
\D_{\uX}^q\uX^2 - q^2\uX^2\D_{\uX}^q = [2]_q\uX + (1-q^2)\sum_{i=1}^m\sum_{\substack{j=1\\ j\neq i}}^m x_j^2e_i\partial_{x_i}^q.
\end{equation}
Axiom (A3) still holds true
\begin{equation}
(\D_{\uX}^q)^2 = - \sum_{i=1}^m (\partial_{x_i}^q)^2 = -\Delta^q
\end{equation}
with the $q$-Laplace operator
\begin{equation}
\Delta^q = \sum_{i_1}^m(\partial_{x_i}^q)^2.
\end{equation}
Using $\uX_q^i = (x_1,...,x_{i-1},qx_i,x_{i+1},...,x_n)$ and $\vert \uX_q^i\vert^2 = x_1^2+...+x_{i-1}^2+q^2x_i^2+x_{i+1}^2+...+x_n^2$ we get
\begin{equation}
\D_{\uX}^q(\uX^2f(\uX)) = [2]_q \sum_{i=1}^m e_ix_if(\uX_q^i) + \vert \uX \vert^2 \D_{\uX}^q(f(\uX))
\end{equation}
or 
\begin{equation}
\D_{\uX}^q(\uX^2f(\uX)) = \sum_{i=1}^m e_i \vert \uX_q^i \vert^2 \partial_{x_i}^q(f(\uX)) + [2]_q\uX f(\uX).
\end{equation}

\subsection{$q$-Euler operator and $q$-Gamma operator}
The $q$-partial derivatives satisfy the Weyl relations
\begin{equation}
\partial_{x_j}^qx_j - qx_j\partial_{x_j}^q = 1.
\end{equation}
Together with the $q$-Dirac operator and the vector variable we can define the $q$-Euler operator $\mathbb{E}^q$. In the continuous case the Euler operator follows from $\D X+X \D = 2\mathbb{E}+m$. Here we get
\begin{equation}
\D_{\uX}^q \uX + \uX \D_{\uX}^q = [2]_q \mathbb{E}^q + m
\end{equation}
with the $q$-Euler operator $\mathbb{E}^q$
\begin{equation}
\mathbb{E}^q = \sum_{i=1}^m x_i\partial_{x_i}^q.
\end{equation}
Now we can study the symmetry relations between $\uX$, $\D_{\uX}^q$, $\mathbb{E}^q$ and $\Delta^q$ using the operations $\{x,y\} = xy+yx$ and $[x,y] = xy-yx$.
\begin{itemize}
\item $\{\uX,\uX\} = -2\vert\uX\vert^2$
\item $\{\D_{\uX}^q,\D_{\uX}^q\} = -2\Delta^q$
\item $\{\D_{\uX}^q,\uX\} = [2]_q\mathbb{E}^q+m$
\item $[\mathbb{E}^q,\uX] = \uX + (q-1)\sum_{i=1}^m x_i^2e_i\partial_{x_i}^q$
\item $[\mathbb{E}_q, \D_{\uX}^q] = \D_{\uX}^q+(q-1) \sum_{i=1}^m x_ie_i(\partial_{x_i}^q)^2$
\item $[\vert \uX \vert^2, \D_{\uX}^q] = [2]_q\uX + (q^2-1) \vert \uX \vert^2 \D_{\uX}^q$
\item $[\mathbb{E}^q,\vert \uX \vert^2] = [2]_q \vert \uX \vert^2 + (q^2-1)\sum_{i=1}^m x_i^3 \partial_{x_i}^q$
\item $[\Delta^q,\uX] = (q^2-1)\sum_{i=1}^m x_ie_i(\partial_{x_i}^q)^2 - [2]_q \D_{\uX}^q$
\item $[\mathbb{E}^q,\Delta^q] = (1-q^2)\sum_{i=1}^m x_i(\partial_{x_i}^q)^3 - [2]_q\Delta^q$
\item $[\Delta^q,\vert \uX \vert^2] =[2]_q[2]_q \mathbb{E}^q + [2]_qm + (q^4-1)\sum_{i=1}^m x_i^2(\partial_{x_i}^q)^2$ 
%(noch falsch, hätte $(\Delta^q)^2$ sein müssen)
\end{itemize}
This results in the usual continuous relations if we take the limit $q\rightarrow 1$.  Similar to the continuous case the $q$-Gamma operator $\Gamma^q$ results out of
\begin{equation}
[\uX,\D_{\uX}^q] = (1-q)\mathbb{E}^q + 2\Gamma^q - m
\end{equation}
with
\begin{equation}
\Gamma^q = - \sum_{i<j} e_ie_j(x_i\partial_{x_j}^q-x_j\partial_{x_i}^q).
\end{equation}
Further the $q$-Gamma operator satisfies 
\begin{equation}
\uX \D_{\uX}^q = \mathbb{E}^q + \Gamma^q.
\end{equation}

Following this we now define a $q$-deformed monogenic homogeneous polynomial as a homogeneous polynomial of degree $k$ $P(\uX) = \vert \uX \vert^k P \left( \frac{\uX}{\vert \uX \vert}\right)$ satisfying $\D_{\uX}^q P = 0$. The space of $q$-deformed monogenic homogeneous polynomials of degree $k$ is then
\begin{equation}
\mathcal{M}^q_k = \{ M(x) \in \mathcal{P}_k: \D_{\uX}^qM = 0\}.
\end{equation}
In the continuous case homogeneous polynomials are the eigenfunctions of the Euler operator $\mathbb{E}P = kP$. The question remains if they are eigenfunctions of the $q$-Euler operator as well. Let $P(\uX)=\sum_{\vert \alpha \vert = k} c_{\alpha}\uX^{\alpha}$ be a homogeneous polynomial of degree $k$ with $\alpha \in \mathbb{N}^m$ a multtindex. Using $P(\uX) = x_1^{\alpha_1}x_2^{\alpha_2}...x_m^{\alpha_m}$ and $|\alpha_1|+...+|\alpha_m| = k$ we get
\begin{align*}
\mathbb{E}^q(P(\uX)) &= \sum_{i=1}^m x_i\partial_{x_i}^q(x_1^{\alpha_1}x_2^{\alpha_2}...x_m^{\alpha_m})\\
&= x_1\frac{(q^{\alpha_1}-1)x_1^{\alpha_1-1}x_2^{\alpha_2}...x_m^{\alpha_m}}{q-1} + x_2\frac{(q^{\alpha_2}-1)x_1^{\alpha_1}x_2^{\alpha_2-1}...x_m^{\alpha_m}}{q-1}\\
&+... + x_m \frac{(q^{\alpha_n}-1)x_1^{\alpha_1}x_2^{\alpha_2}...x_m^{\alpha_m-1}}{q-1}\\
&= [\alpha_1]_qx_1^{\alpha_1}x_2^{\alpha_2}...x_m^{\alpha_m}+ [\alpha_2]_q x_1^{\alpha_1}x_2^{\alpha_2}...x_m^{\alpha_m}+...+[\alpha_m]_q x_1^{\alpha_1}x_2^{\alpha_2}...x_m^{\alpha_m}\\
&= ([\alpha_1]_q+[\alpha_2]_q+...+[\alpha_m]_q)P(\uX).
\end{align*}
Therefore homogeneous polynomials are eigenfunctions of the $q$-Euler operator to the eigenvalue $[\alpha_1]_q+ \ldots +[\alpha_m]_q$. As $\alpha_1+\ldots+\alpha_m = k$ this results in $\mathbb{E}(P(\uX) = kP(\uX)$ for $q \rightarrow 1$.
Depending on the order $k$ and the specific partition of $k$, we obtain not only the eigenvalue $k = [1]_q + \ldots + [1]_q$ because $[l]_q = 1 + q+ \ldots +q^{l-1}, l\in \mathbb{N}.$

\begin{example}  For $k=3$ and $m=3$ we have the following partitions which leads to the corresponding eigenvalues. \\[1ex]
\begin{tabular}{l|l|l}
Homogeneous polynomials  & Partition & Eigenvalue \\ \hline
$x_1x_2x_3$  & $3 = 1+1+1$   &    $[1]_q+[1]_q+[1]_q= 1+1+1 = 3$ \\
$x_1x_2^2,  x_1x_3^2, x_2x_3^2,  x_1^2x_2, x_1^2x_3, x_2^2x_3 $ & $3=1+2 $    &      $ [1]_q+[2]_q=1+ 1+q = 2+q $ \\
$x_1^3, x_2^3, x_3^3$ & $3 $ &                      $[3]_q= 1+q+q^2$ \\
\end{tabular}
\end{example}

\section{Fischer decomposition}
Using multi-index notation $\alpha = (\alpha_1,\ldots,\alpha_m) \in \mathbb{N}^m$ we have
\begin{itemize}
\item $\uX^{\alpha} = x_1^{\alpha_1}\cdots x_m^{\alpha_m}$
\item $\alpha! = \alpha_1! \cdots \alpha_n!$
\item $\vert \alpha \vert = \sum_{i=1}^m \alpha_i$
\item $(\partial_{\uX}^q)^{\alpha} = (\partial_{x_1}^q)^{\alpha_1}\cdots (\partial_{x_m}^q)^{\alpha_m}$.
\end{itemize}
We use the standard basis for $k$-homogeneous Clifford-valued polynomials $\mathcal{P}_k = \{\uX^{\alpha}: \vert \alpha \vert = k\}$. First we define an inner product on the complex vector space $\mathcal{P}_k$ as follows.\\
For $R_1,\,R_2 \in \mathcal{P}_k$ with $R_i(\uX) = \sum_{\vert \alpha \vert =k} x^{\alpha}a_{\alpha}^i,\, a_{\alpha}^i \in \mathcal{C}\ell_{0,m},\, i=1,2$ we define
\begin{equation}
\langle R_1,R_2\rangle_{k,q} = \sum_{\vert \alpha \vert = k} [\alpha]_q!(\overline{a_{\alpha}^1}a_{\alpha}^2)_0.
\end{equation}
This is a scalar product, the so called Fischer inner product.
\begin{satz}
For $R_1,\,R_2 \in \mathcal{P}_k$ we get
\begin{equation}
\langle R_1, R_2\rangle_{k,q} = (\overline{R}_1(\partial_{\uX}^q)R_2)_0
\label{eq:fischer2}
\end{equation}
where $R_1(\partial_{\uX}^q)$ denotes the operator obtained by replacing the $x_j$ in $R_1$ with the $q$-partial derivative $\partial_{x_j}^q$.
\end{satz}
Proof: We use the fact that
\begin{align*}
(\partial_{\uX}^q)^{\alpha}\uX^{\beta} = \begin{cases}
	[\alpha]_q!, & \alpha = \beta \\
	0	,	    & \alpha \neq \beta.
\end{cases}
\end{align*}
Using $R_1, R_2 \in \mathcal{P}_k$ with coefficients $a_{\alpha}^i = 1$ we get for the components in (\ref{eq:fischer2})
\begin{align}
R_1(\partial_{\uX}^q) &= (\partial_{x_1}^q)^{\alpha_1}...(\partial_{x_n}^q)^{\alpha_n},\text{ } \alpha_1+...+\alpha_n = k\\
R_2(\uX) &= x_1^{\beta_1}...x_n^{\beta_n},\text{ } \beta_1+...+\beta_n = k.
\end{align}
There are three different cases to examine. It is sufficient to verify these with a single variable $x_i$.\\
The first case is $\alpha_i = \beta_i$ for each $i=1,...,n$
\begin{align*}
(\partial_{x_i}^q)^{\alpha_i}x_i^{\alpha_i} &= (\partial_{x_i}^q)^{\alpha_i-1}\left(\frac{q^{\alpha_i}x_i^{\alpha_i}-x_i^{\alpha_i}}{(q-1)x_i}\right)\\
&=(\partial_{x_i}^q)^{\alpha_i-1}([\alpha_i]_qx_i^{\alpha_i-1})\\
&= \ldots \\
&= [\alpha_i]_q!.
\end{align*}
The next case is $\alpha_i > \beta_i$
\begin{align*}
(\partial_{x_i}^q)^{\alpha_i+1}x_i^{\alpha_i} &= \partial_{x_i}^q([\alpha_i]_q!)\\
&= \frac{[\alpha_1]_q!-[\alpha_i]_q!}{(q-1)x_i}\\
&= 0.
\end{align*}
Finally we have $\alpha_i<\beta_i$
\begin{align*}
(\partial_{x_i}^q)^{\alpha_i} x_i^{\beta_i} \vert_{x_i=0} &= [\alpha_i]_q!x_i^{\beta_i-\alpha_i} \vert_{x_i=0}\\
&= 0.
\end{align*}
Together we get for $R_1,R_2 \in \mathcal{P}_k$ with any coefficients $a_{\alpha}^i\in \mathcal{C}\ell_{0,n}$
\begin{align*}
(\overline{R_1}(\partial_{\uX}^q)R_2)_0 &= \sum_{\vert\alpha\vert = k} [\alpha]_q!(\overline{a_{\alpha}^1}a_{\alpha}^2)_0\\
&= \langle R_1,R_2\rangle_{k,q} 
\end{align*}
\eop\\
\begin{satz}
For all $Q \in \mathcal{P}_k$ and $P \in \mathcal{P}_{k+1}$
\begin{equation}
\langle \uX Q,P\rangle_{k+1,q} = - \langle Q, \D_{\uX}^qP\rangle_{k,q}.
\end{equation}
\end{satz}
Proof: Due to linearity, it is sufficient only to consider monomials
\begin{align}
Q(\uX) &= e_B x_1^{\beta_1}...x_n^{\beta_n}=x_1^{\beta_1}...x_n^{\beta_n}e_B,\text{ } \beta_1+...+\beta_n = k\\
P(\uX) &= e_A x_1^{\alpha_1}...x_n^{\alpha_n}=x_1^{\alpha_1}...x_n^{\alpha_n}e_A,\text{ } \alpha_1+...+\alpha_n = k+1
\end{align}
Further we only need to consider the variable $x_i$. On the left hand side we get
\begin{align*}
\langle e_ix_iQ,P\rangle_{k+1,q} &= (\overline{e_i\partial_{x_i}^q(\partial_{x_1}^q)^{\beta_1}...(\partial_{x_i}^q)^{\beta_i}...(\partial_{x_n}^q)^{\beta_n}e_B}x_1^{\alpha_1}...x_n^{\alpha_n}e_A)_0\\
&= (\overline{e_i(\partial_{x_1}^q)^{\beta_1}...(\partial_{x_i}^q)^{\beta_i+1}...(\partial_{x_n}^q)^{\beta_n}e_B}x_1^{\alpha_1}...x_n^{\alpha_n}e_A)_0\\
&=(\overline{(\partial_{x_1}^q)^{\beta_1}...(\partial_{x_i}^q)^{\beta_i+1}...(\partial_{x_n}^q)^{\beta_n}e_B}\overline{e_i}x_1^{\alpha_1}...x_n^{\alpha_n}e_A)_0.
\end{align*}
While the right hand side is
\begin{align*}
\langle Q, e_i\partial_{x_i}^qP\rangle_{k,q} &= (\overline{(\partial_{x_1}^q)^{\beta_1}...(\partial_{x_n}^q)^{\beta_n}e_B}e_i\partial_{x_i}^qx_1^{\alpha_1}...x_i^{\alpha_i}...x_n^{\alpha_n}e_A)_0\\
&= (\overline{(\partial_{x_1}^q)^{\beta_1}...(\partial_{x_n}^q)^{\beta_n}e_B} e_i [\alpha_i]_q x_1^{\alpha_1}...x_i^{\alpha_i-1}...x_n^{\alpha_n}e_A)_0.
\end{align*}
Because $\overline{e_i} = -e_i$ we get either
\begin{equation*}
\langle \uX Q,P\rangle_{k+1,q} = -\langle Q,\D_{\uX}^qP\rangle_{k,q}
\end{equation*}
or both sides equal zero.\eop\\
\\
This allows the following theorem.
\begin{satz}
For $k \in \mathbb{N}$ we have
\begin{equation}
\mathcal{P}_k = \mathcal{M}^q_k \oplus \uX\mathcal{P}_{k-1}.
\end{equation}
Further the subspaces $\mathcal{M}^q_k$ and $\uX\mathcal{P}_{k-1}$ of $\mathcal{P}_k$ are orthogonal with respect ot the scalar product $\langle \cdot,\cdot \rangle_{k,q}$.
\end{satz}
Proof: As $\mathcal{P}_k = \uX\mathcal{P}_{k-1} \oplus (\uX\mathcal{P}_{k-1})^{\bot}$ it suffices to proof $\mathcal{M}^q_k = (\uX\mathcal{P}_{k-1})^{\bot}$. For the first inclusion take any $R_{k-1} \in \mathcal{P}_{k-1}$ and $R_k \in \mathcal{P}_k$. Suppose now that
\begin{equation*}
\langle \underbrace{\uX R_{k-1}}_{\in \uX\mathcal{P}_{k-1}},\underbrace{R_k}_{\in(\uX\mathcal{P}_{k-1})^{\bot}} \rangle_{k,q}= 0.
\end{equation*}
Then we have $\langle R_{k-1},\D_{\uX}^qR_k \rangle_{k_1,q} = 0$ for each $R_{k-1} \in \mathcal{P}_{k-1}$ (see theorem 2). We put $R_{k-1} = \D_{\uX}^qR_k$ and it follows that $\D_{\uX}^qR_k = 0$. Therefore $R_k \in \mathcal{M}_k$ and finally $(\uX\mathcal{P}_{k-1})^{\bot} \subset \mathcal{M}^q_k$.\\
For the other inclusion take $P_k \in \mathcal{M}^q_k$. For each $R_{k-1} \in \mathcal{P}_{k-1}$
\begin{align*}
\langle \uX R_{k-1},P_k\rangle_{k,q} &= -\langle R_{k-1},\underbrace{\D_{\uX}^qP_k}_{=0}\rangle_{k-1,q}\\
&= 0.
\end{align*}
As $\uX R_{k-1} \in \uX\mathcal{P}_{k-1}$ we have $P_k \in (\uX\mathcal{P}_{k-1})^{\bot}$ and accordingly $\mathcal{M}^q_k \subset (\uX\mathcal{P}_{k-1})^{\bot}$. Therefore $\mathcal{M}^q_k = (\uX\mathcal{P}_{k-1})^{\bot}$. \eop\\
\\
Theorem 3 allows a Fischer decomposition of the space of homogeneous polynomials $\mathcal{P}_k$:
\begin{equation}
\mathcal{P}_k = \sum_{s=0}^k \oplus \uX^s\mathcal{M}^q_{k-s}.
\end{equation}
This follows directly from theorem 3 using
\begin{equation*}
\uX\mathcal{P}_{k-1} = \uX\mathcal{M}^q_{k-1} \oplus \uX^2\mathcal{P}_{k-2},\, \uX^2\mathcal{P}_{k-2} = \uX^2\mathcal{M}^q_{k-2} \oplus \uX^3\mathcal{P}_{k-3},...\, .
\end{equation*}

\section{Cauchy-Kovalevskaya extension}
To generate monogenic functions, we can use the Cauchy-Kovalevskaya extension theorem.
\begin{satz}
Let $f(x_1,...,x_m) = f(\uX)$ in $\mathbb{R}^m$. The $q$-deformed Cauchy-Kovalevskaya extension of the function $f(\uX)$ is
\begin{equation}
\text{CK}(f(\uX)) = f^*(x_0,x_1,...,x_m) = \sum_{k=0}^{\infty} \frac{1}{[k]_q!}x_0^k(\overline{e}_0\D_{\uX}^q)^kf(\uX).
\end{equation}
$f^*$ is a monogenic function satisfying $\D_x^q(f^*(x_0,...,x_m)) = 0$, with $\D_x^q = -\sum_{i=0}^m e_i\partial_{x_i}^q = -e_0\partial_{x_0}^q + \D_{\uX}^q$ and $e_0^2 = -1$. Further $f^*\vert_{x_0=0} = f$.
\end{satz}

\begin{example}[$q$-deformed Fueter variables]
Applying the CK-extension to the function $x_l$ in $\mathbb{R}^m$ we get
\begin{align*}
f^*(x_0,x_1,...,x_l,...,x_m) &= \sum_{k=0}^{\infty} \frac{1}{[k]_q!}x_0^k(\overline{e}_0\D_{\uX}^q)^k(x_l)\\
&=\frac{1}{[0]_q!}x_0^0(\overline{e}_0\D_{\uX}^q)^0(x_l) + \frac{1}{[1]_q!}x_0(\overline{e}_0\D_{\uX}^q)(x_l)+0+...\\
&= x_l - x_0\overline{e}_0e_l.
\end{align*}
For $x_0 = 0$ the CK-extension $f^*$ reduces to the original function $f$. Clearly this is a monogenic function:
\begin{align*}
\D_x^q(x_l-x_0\overline{e}_0e_l) &= -e_0\partial_{x_0}^q(x_l-x_0\overline{e}_0e_l) - e_l\partial_{x_l}^q(x_l-x_0\overline{e}_0e_l)\\
&= e_0\overline{e}_0e_l - e_l = 0.
\end{align*}
\end{example}

\begin{example}
We want to compute the Cauchy-Kovalevskaya extension of $x_ix_j$. For $i\not= j$ we obtain
\begin{align*}
(x_ix_j)^* &= \sum_{k=0}^{\infty} \frac{1}{[k]_q!}x_0^k(\overline{e}_0\D_{\uX}^q)^k(x_ix_j)\\
&=\frac{1}{[0]_q!}x_0^0(\overline{e}_0\D_{\uX}^q)^0(x_ix_j) + \frac{1}{[1]_q!}x_0(\overline{e}_0\D_{\uX}^q)(x_ix_j)+\frac{1}{[2]_q!}x^2_0(\overline{e}_0\D_{\uX}^q)^2(x_ix_j) + \ldots  \\
& = x_ix_j - x_0\overline{e}_0(x_je_i+ x_ie_j) + \frac{1}{[2]_q!}x^2_0 (\overline{e}_0(-\overline{e}_0))(e_je_i + e_ie_j) \\
& =  x_ix_j - x_0\overline{e}_0(x_je_i+ x_ie_j) \\
&= \frac{1}{2!}\left( (x_i - x_0\overline{e}_0e_i)(x_j - x_0\overline{e}_0e_j) + (x_j - x_0\overline{e}_0e_j)(x_i - x_0\overline{e}_0e_i) \right).
\end{align*}
But in case of $i=j$ we get
\begin{align*}
(x_i^2)^* &= \sum_{k=0}^{\infty} \frac{1}{[k]_q!}x_0^k(\overline{e}_0\D_{\uX}^q)^k(x_i^2)\\
&=\frac{1}{[0]_q!}x_0^0(\overline{e}_0\D_{\uX}^q)^0(x_i^2) + \frac{1}{[1]_q!}x_0(\overline{e}_0\D_{\uX}^q)(x_i^2)+\frac{1}{[2]_q!}x^2_0(\overline{e}_0\D_{\uX}^q)^2(x_i^2) + \ldots  \\
& = x_i^2 - x_0\overline{e}_0( [2]_q x_ie_i) +\frac{1}{[2]_q!}x^2_0 (\overline{e}_0(-\overline{e}_0))([2]_qe_i^2) \\
& =  x_i^2 - [2]_qx_0x_i\overline{e}_0e_i -x_0^2 \not = (x_i-x_0\overline{e}_0e_i)^2,\quad \text{because  } [2]_q = 1+q \not= 2 \text{ for } q\not= 1.
\end{align*}
\end{example}
Therefore, the Cauchy-Kovalevskaya extension of products is not in general the product of $q$-deformed Fueter variables. 
\printbibliography

\end{document}